\documentclass[12pt]{article}                                               
\input{amssym.def}                                                           
\input{amssym.tex}                                                           
                                                                             
\begin{document}                                                             
%**********************************************                              
\title{$AF$-algebras and  topology of mapping tori}

\author{Igor  ~Nikolaev
\footnote{Partially supported 
by NSERC.}}

%**************************************************

\date{}
 \maketitle

%**************************************************

\newtheorem{thm}{Theorem}
\newtheorem{lem}{Lemma}
\newtheorem{dfn}{Definition}
\newtheorem{rmk}{Remark}
\newtheorem{cor}{Corollary}
\newtheorem{cnj}{Conjecture}
\newtheorem{exm}{Example}

%*************************************************

\centerline{{\it  In memory of  ~D. ~V. ~Anosov}}

%******************************************************************
\begin{abstract}
A covariant functor from the category of mapping tori 
to a category of $AF$-algebras is constructed;  the functor takes continuous maps
between such manifolds   to stable homomorphisms between  the corresponding  
$AF$-algebras.  We use this  functor to develop an obstruction theory for the 
torus bundles of dimension $2, 3$ and $4$.

\vspace{7mm}

{\it Key words and phrases:  Anosov diffeomorphism,  $AF$-algebra}

\vspace{5mm}
{\it MSC:   46L85 (noncommutative topology);  55S35 (obstruction theory)}
\end{abstract}

%**************************************************************************
\section*{Introduction}
%***************************************************************************
This paper studies applications of operator algebras in topology;
the operator algebras in question are the so-called {\it $AF$-algebras}
and the topological spaces are certain {\it mapping tori},  i.e.  circle bundles
with a fiber $M$ and monodromy  $\varphi: M\to M$.  Recall that
a very fruitful  approach to topology consists in construction of maps (functors) from the 
topological spaces to certain algebraic objects,  so that continuous maps
between the spaces become homomorphisms of the corresponding algebraic 
entities.  The functors usually take value in the finitely generated groups 
(abelian or not) and, therefore, reduce topology  to a simpler algebraic problem. 
The rings of operators on a Hilbert space are neither finitely generated nor commutative 
and, at the first glance, if ever such a reduction exists,  it will not simplify the problem.  
Yet  it is not so:   we define an operator algebra, the so-called fundamental 
$AF$-algebra, which  yields a set of simple obstructions (invariants)
to existence of continuous maps in a  class of  manifolds fibering over  the
circle.  One obstruction turns out to be the Galois group of the fundamental
$AF$-algebra; this invariant dramatically simplifies for a class of the so-called
tight torus bundles, so that topology boils down to a division test for a finite
set of natural numbers.

\medskip\noindent
{\bf A. $AF$-algebras ([Bratteli 1972]  \cite{Bra1}).}
The $C^*$-algebra $A$ is an algebra over the complex numbers $\Bbb C$ with a norm 
$a\mapsto ||a||$ and an involution $a\mapsto a^*, a\in A$, such that $A$ is
complete with respect to the norm, and such that $||ab||\le ||a||~||b||$ 
and $||a^*a||=||a||^2$ for every $a,b\in A$. 
Any commutative $C^*$-algebra $A$  is isomorphic to the $C^*$-algebra 
$C_0(X)$ of continuous complex-valued functions on a locally compact Hausdorff space $X$
vanishing at the infinity; 
the algebras which are not commutative are deemed as   noncommutative topological
spaces. A {\it stable homomorphism} $A\to A'$ is defined as the (usual) homomorphism
$A\otimes {\cal K}\to A'\otimes {\cal K}$, where ${\cal K}$ is the $C^*$-algebra
of compact operators on a Hilbert space; such a homomorphism corresponds to
 a continuous map between the noncommutative spaces $A$ and $A'$.    
 The matrix algebra $M_n({\Bbb C})$ is an example of noncommutative finite-dimensional
$C^*$-algebra;  a natural generalization  are  approximately finite-dimensional ($AF$-) algebras, 
which  are given by an ascending sequence  
$M_1\buildrel\rm\varphi_1\over\longrightarrow M_2\buildrel\rm\varphi_2\over\longrightarrow\dots$
of finite-dimensional semi-simple $C^*$-algebras $M_i=M_{n_1}({\Bbb C})\oplus\dots\oplus M_{n_k}({\Bbb C})$
and homomorphisms $\varphi_i$  arranged into an infinite graph 
 as follows.  The two sets of vertices $V_{i_1},\dots, V_{i_k}$ and $V_{i_1'},\dots, V_{i_k'}$
are joined by the $b_{rs}$ edges, whenever the summand $M_{i_r}$ contains $b_{rs}$
copies of the summand $M_{i_s'}$ under the embedding $\varphi_i$; as $i\to\infty$,
one gets a {\it Bratteli diagram} of the $AF$-algebra. Such a diagram is defined by
an infinite sequence of {\it incidence matrices}  $B_i=(b_{rs}^{(i)})$.  If the homomorphisms 
$\varphi_1 =\varphi_2=\dots=Const$,  the $AF$-algebra is called {\it stationary}; 
its Bratteli diagram looks like an infinite graph  with the incidence matrix 
$B=(b_{rs})$ repeated over and over again.

\medskip\noindent
{\bf B. $AF$-algebra of measured foliation.}
Let $M$ be a compact manifold of dimension $m$ and
${\cal F}$ a codimension $k$ measured foliation of $M$;  it is known that ${\cal F}$
is tangent to the hyperplane $\omega(p)=0$
at each point $p\in M$, where $\omega\in H^k(M; {\Bbb R})$ is a closed $k$-form, 
see e.g. [Plante 1975]  \cite{Pla1}. Denote by  $\lambda_i>0$  the  periods of $\omega$
against a basis in the  homology group $H_k(M)$ and   
consider the vector $\theta=(\theta_1,\dots,\theta_{n-1})$,
where $\theta_i=\lambda_{i+1}/\lambda_1$ and $n=rank~H_k(M)$. 
Let $\lim_{i\to\infty}B_i$ be the
Jacobi-Perron continued fraction convergent to the vector $(1,\theta)$;
here $B_i\in GL_n({\Bbb Z})$ are the non-negative matrices with $det~(B_i)=1$ [Bernstein 1971]  \cite{B}.
An $AF$-algebra ${\Bbb A}_{\cal F}$ is called {\it associated} to ${\cal F}$,
if its  Bratteli  diagram is given by the matrices $B_i$;
the Bratteli diagram defines an isomorphism class of ${\Bbb A}_{\cal F}$ [Bratteli 1972]  \cite{Bra1}. 
The algebra ${\Bbb A}_{\cal F}$ has a spate of remarkable 
properties, e.g.  the topologically conjugate (or, induced) foliations 
have stably isomorphic (or, stably homomorphic) $AF$-algebras
(lemma \ref{lm1});  the  dimension group of 
${\Bbb A}_{\cal F}$ [Effros 1981]  \cite{E} coincides with the Plante group 
$P({\cal F})$ of foliation ${\cal F}$ [Plante 1975]  \cite{Pla1}.

\medskip\noindent
{\bf C. Fundamental $AF$-algebras and main result.}
Let $\varphi: M\to M$ be an Anosov diffeomorphism of $M$ [Anosov 1967]  \cite{Ano1}; if $p$ is a fixed point
of $\varphi$, then $\varphi$ defines an invariant measured foliation ${\cal F}$ of $M$ given by the stable manifold 
$W^s(p)$ of $\varphi$ at the point  $p$ [Smale 1967]  \cite{Sma1}, p.760.  The associated $AF$-algebra ${\Bbb A}_{\cal F}$ is 
stationary (lemma \ref{lm2});  we call the latter a {\it fundamental $AF$-algebra} and denote it by
${\Bbb A}_{\varphi}:={\Bbb A}_{\cal F}$.   Consider the mapping torus of $\varphi$, i.e. a manifold  
$M_{\varphi}:= M\times [0,1]~/~\sim$, where $(x,0)\sim (\varphi(x), 1), ~\forall x\in M$.  
Let ${\cal M}$ be a category of the mapping tori
of all Anosov's  diffeomorphisms; 
the arrows  of ${\cal M}$ are continuous maps between the mapping tori. Likewise, let ${\cal A}$
be a category of all fundamental $AF$-algebras; the arrows  of ${\cal A}$
are stable homomorphisms between the fundamental $AF$-algebras. By $F: {\cal M}\to {\cal A}$
we understand a map given by the formula $M_{\varphi}\mapsto {\Bbb A}_{\varphi}$, where $M_{\varphi}\in {\cal M}$
and ${\Bbb A}_{\varphi}\in {\cal A}$. Our main result can be stated as follows.
%*************************************************************************
\begin{thm}\label{thm1}
The map $F$ is a functor, which sends each continuous map $N_{\psi}\to M_{\varphi}$ to a stable
homomorphism ${\Bbb A}_{\psi}\to {\Bbb A}_{\varphi}$ of the corresponding fundamental $AF$-algebras.
\end{thm}
%***********************************************************************************

\medskip\noindent
{\bf D. Applications.}
Theorem \ref{thm1} has a natural application,
since  stable homomorphisms of the fundamental $AF$-algebras are easier 
to detect,  than  continuous maps between  manifolds $N_{\psi}$ and $M_{\varphi}$;
such homomorphisms are in bijection  with the inclusions of certain ${\Bbb Z}$-modules
lying  in a  (real) algebraic number field.  Often it is possible to prove, that no inclusion
is possible and, thus, draw a  topological conclusion about the  maps (an obstruction
theory). Namely,  since ${\Bbb A}_{\psi}$ is stationary, it has a constant incidence matrix $B$;
the splitting field of the polynomial $det~(B-xI)$ we denote by $K_{\psi}$ and call $Gal~(K_{\psi}|{\Bbb Q})$ the {\it Galois group} of the algebra ${\Bbb A}_{\psi}$.  
 Suppose that $h: {\Bbb A}_{\psi}\to {\Bbb A}_{\varphi}$ is a stable homomorphism;
since the corresponding invariant foliations ${\cal F}_{\psi}$ and ${\cal F}_{\varphi}$ are induced,
their Plante groups are included $P({\cal F}_{\varphi})\subseteq P({\cal F}_{\psi})$ 
and, therefore,  ${\Bbb Q}(\lambda_{B'})\subseteq K_{\psi}$,
where $\lambda_{B'}$ is the Perron-Frobenius eigenvalue of the matrix $B'$ attached to
${\Bbb A}_{\varphi}$.  Thus,  stable homomorphisms are in bijection   with  subfields of the algebraic
number field $K_{\psi}$;  their  classification achieves perfection in terms of the Galois theory,
since the subfields  are in a  one-to-one correspondence  with the subgroups of  
$Gal~({\Bbb A}_{\psi})$ [Morandi 1996]  \cite{M}.   
In particular, when $Gal~({\Bbb A}_{\psi})$ is simple, there are only trivial stable
homomorphisms; thus, the structure of  $Gal~({\Bbb A}_{\psi})$ is  an obstruction  (an invariant)
to existence of a continuous map between the manifolds $N_{\psi}$ and $M_{\varphi}$.
Is our invariant effective? The answer is positive for a class of the so-called tight
torus bundles; in this case $N_{\psi}$ is given by a monodromy matrix,  which  is similar to the
matrix $B$.  The obstruction theory for the tight torus bundles of any dimension can 
be  completely determined;  it reduces to  a  divisibility test for  a finite set of 
natural numbers.  For the sake of clarity, the test is done in dimension $m=2,3$ and $4$ and followed 
by the numerical examples.  
%******************************************************************************
\begin{rmk}
\textnormal{
Notice that for the tight torus bundles (see Section 3.2) our results
can be proved using the theory of hyperbolic diffeomorphism $\psi: T^m\to T^m$
alone; however, our approach seems to be more general leading to the essentially
new topological invariants. 
}
\end{rmk}
%********************************************************************************  

%***********************************************************
\tableofcontents
%*********************************************************** 

%**************************************************************************
\section{Preliminaries}
%***************************************************************************
%**************************************************************************
\subsection{Measured foliations}
%***************************************************************************
By a $q$-dimensional,  class  $C^r$ foliation of an $m$-dimensional manifold $M$ 
one understands a decomposition of $M$ into a union of
disjoint connected subsets $\{ {\cal L}_{\alpha}\}_{\alpha\in A}$, called
the {\it leaves} [Lawson 1974]  \cite{Law1}.  They must satisfy the
following property: each point in $M$ has a neighborhood $U$
and a system of local class $C^r$ coordinates 
$x=(x^1,\dots, x^m): U\to {\Bbb R}^m$ such that for each leaf 
${\cal L}_{\alpha}$, the components of $U\cap {\cal L}_{\alpha}$
are described by the equations $x^{q+1}=Const, \dots, x^m=Const$.
Such a foliation is denoted by ${\cal F}=\{ {\cal L_{\alpha}}\}_{\alpha\in A}$.
The number $k=m-q$ is called the  codimension of the foliation.
An  example of a codimension $k$ foliation ${\cal F}$ 
is given by a closed $k$-form $\omega$ on $M$: the leaves of ${\cal F}$
are tangent to the hyperplane  $\omega(p)=0$ at each point $p$ of $M$. 
The $C^r$-foliations ${\cal F}_0$ and ${\cal F}_1$ 
of codimension $k$
are said to be $C^s$-conjugate ($0\le s\le r$), if there exists an (orientation-preserving) 
diffeomorphism of $M$, of class $C^s$, which maps the leaves of ${\cal F}_0$
onto the leaves of ${\cal F}_1$;  when $s=0$, ${\cal F}_0$ and ${\cal F}_1$
are {\it topologically conjugate}. 
Denote by  $f: N\to M$  a map of class $C^s$ ($1\le s\le r$) of a manifold
$N$ into $M$;  the map $f$ is said to be {\it transverse} to ${\cal F}$,
if for all $x\in N$ it holds $T_y(M)=\tau_y({\cal F})+f_* T_x(N)$,
where $\tau_y({\cal F})$ are the vectors of $T_y(M)$ tangent
to ${\cal F}$ and  $f_*: T_x(N)\to T_y(M)$ is  the linear map
on tangent vectors induced by $f$, where $y=f(x)$. 
If the map  $f: N\to M$ is transverse to a foliation ${\cal F}'=\{{\cal L}\}_{\alpha\in A}$
on $M$, then $f$  induces a class $C^s$ foliation ${\cal F}$ on $N$,
where the leaves are defined as $f^{-1}({\cal L}_{\alpha})$ for all 
$\alpha\in A$; it is immediate, that $codim~({\cal F})=codim~({\cal F}')$. 
We shall call ${\cal F}$ an {\it induced foliation}.  
When $f$ is a submersion, it is transverse  to any foliation of $M$;
in this case, the induced foliation ${\cal F}$ is correctly defined for all ${\cal F}'$
on $M$ [Lawson 1974]  \cite{Law1}, p.373. 
Notice,  that for $M=N$  the above definition corresponds to  
topologically conjugate foliations ${\cal F}$ and ${\cal F}'$.  
To introduce measured foliations, denote by $P$ and $Q$  two $k$-dimensional
submanifolds of $M$,  which are everywhere transverse to a foliation ${\cal F}$ of codimension $k$. 
Consider a collection of $C^r$ homeomorphisms between subsets of $P$ and $Q$
induced by  a return map along the leaves of ${\cal F}$.
The collection of all such homeomorphisms between subsets of all possible pairs of 
transverse manifolds generates a {\it holonomy pseudogroup} of ${\cal F}$
under composition of the homeomorphisms [Plante 1975]  \cite{Pla1}, p.329.  
A foliation ${\cal F}$  is said to have measure preserving holonomy,
if its holonomy pseudogroup has a non-trivial invariant measure, which is
finite on compact sets; for brevity, we call ${\cal F}$ a  {\it measured foliation}.
An example of measured foliation is a foliation, determined by 
the closed $k$-form $\omega$; the restriction of $\omega$ to a transverse $k$-dimensional
manifold determines a volume element, which gives a positive invariant measure on open sets.     
Each measured foliation ${\cal F}$ defines an element of the cohomology
group $H^k(M; {\Bbb R})$ [Plante 1975]  \cite{Pla1};  in the case of ${\cal F}$ given by a closed
$k$-form $\omega$, such an element coincides with the de Rham cohomology class of $\omega$,  
{\it ibid}. In view of the isomorphism $H^k(M; {\Bbb R})\cong Hom~(H_k(M), {\Bbb R})$,
foliation ${\cal F}$ defines a linear map $h$ from the $k$-th homology group $H_k(M)$ 
to ${\Bbb R}$;  by the {\it Plante group} $P({\cal F})$
we shall understand a finitely generated abelian subgroup $h(H_k(M)/Tors)$ of the real line
${\Bbb R}$.  If $\{\gamma_i\}$ is a basis of the homology group $H_k(M)$,
then the periods $\lambda_i=\int_{\gamma_i}\omega$ are  generators of the group  $P({\cal F})$ [Plante 1975]  \cite{Pla1}.

%**************************************************************************
\subsection{$AF$-algebra of measured foliation}
%***************************************************************************
Let  $\lambda=(\lambda_1,\dots,\lambda_n)$ be a basis of the Plante group $P({\cal F})$
of  a  measured foliation ${\cal F}$, such that $\lambda_i>0$.  
Take a vector $\theta=(\theta_1,\dots,\theta_{n-1})$
with $\theta_i=\lambda_{i+1} / \lambda_1$; the Jacobi-Perron continued fraction of
vector $(1, \theta)$ (or, projective class of vector $\lambda$)  is given by the formula ([Bernstein 1971]  \cite{B}, p.13):
%**************************************************************
\begin{equation}\label{eq1}
\left(\matrix{1\cr \theta}\right)=
\lim_{i\to\infty} \left(\matrix{0 & 1\cr I & b_1}\right)\dots
\left(\matrix{0 & 1\cr I & b_i}\right)
\left(\matrix{0\cr {\Bbb I}}\right)=
\lim_{i\to\infty} B_i\left(\matrix{0\cr {\Bbb I}}\right),
\end{equation}
%****************************************************************
where $b_i=(b^{(i)}_1,\dots, b^{(i)}_{n-1})^T$ is a vector of  non-negative integers,  
$I$ the unit matrix and ${\Bbb I}=(0,\dots, 0, 1)^T$;  the $b_i$ are obtained from $\theta$
by the Euclidean algorithm, see [Bernstein 1971]  \cite{B}, pp.2-3 for details.
An $AF$-algebra given by the Bratteli diagram with the incidence matrices $B_i$
(and unital homomorphisms $M_i\to M_{i+1}$)
will be called {\it associated} to the  foliation ${\cal F}$; we shall denote such 
an algebra by ${\Bbb A}_{\cal F}$.  
Taking another basis of the Plante group $P({\cal F})$ gives 
an $AF$-algebra  which is stably isomorphic to  ${\Bbb A}_{\cal F}$;  
this is an algebraic recast of the main  property  of the Jacobi-Perron fractions. 
It is known,  that the Bratteli diagram defines the $AF$-algebra up to an
isomorphism [Bratteli 1972]  \cite{Bra1};  by ${\Bbb A}_{\cal F}$  we mean the isomorphism class.  
Note, that if ${\cal F}'$ is a measured foliation on a manifold $M$ and 
$f:N\to M$ is a submersion, the induced foliation 
${\cal F}$ on $N$ is a measured foliation.
We shall denote by ${\cal MF}ol$ the category of all manifolds
with measured foliations (of fixed codimension),  whose arrows are submersions of
the manifolds;  by ${\cal M}_0{\cal F}ol$ we understand  a subcategory of ${\cal MF}ol$, 
consisting  of manifolds, whose foliations have a  unique transverse measure. 
Let   ${\cal R}ng$ be the  category 
of the (isomorphism classes of) $AF$-algebras given by {\it convergent} Jacobi-Perron fractions (\ref{eq1}),
so that the arrows of ${\cal R}ng$  are the stable homomorphisms of the $AF$-algebras.
By $F$ we denote a map between  ${\cal M}_0{\cal F}ol$ and ${\cal R}ng$ given by the
formula ${\cal F}\mapsto {\Bbb A}_{\cal F}$.  Notice, that $F$ is correctly defined,
since foliations with the unique measure have the convergent Jacobi-Perron
fractions; this assertion follows from [Bauer 1996]  \cite{Bau1}. 
%************************************************************************
\begin{lem}\label{lm1}
The map  $F: {\cal M}_0{\cal F}ol\longrightarrow  {\cal R}ng$ is a functor,
which sends any pair of induced foliations to a pair of stably homomorphic
$AF$-algebras.
\end{lem}
%***********************************************************************

%**************************************************************************
\subsection{Fundamental $AF$-algebras}
%***************************************************************************
Let $M$ be an $m$-dimensional manifold and $\varphi: M\to M$ a diffeomorphisms  of $M$;
recall, that  an orbit of point $x\in M$ is the subset $\{\varphi^n(x) ~|~n\in {\Bbb Z}\}$
of $M$. The finite orbits $\varphi^m(x)=x$ are called periodic; 
when $m=1$, $x$ is a {\it fixed point} of diffeomorphism $\varphi$.
The fixed point $p$ is {\it hyperbolic} if the eigenvalues $\lambda_i$ of 
the linear map $D\varphi(p): T_p(M)\to T_p(M)$ do not lie on the unit circle. 
If $p\in M$ is a hyperbolic fixed point of a diffeomorphism $\varphi: M\to M$,
denote by  $T_p(M)=V^s+V^u$ the corresponding decomposition of the tangent space 
under the linear map $D\varphi(p)$, where $V^s$ ($V^u$) is the eigenspace of $D\varphi(p)$
corresponding to $|\lambda_i|>1$ ($|\lambda_i|<1$). 
For a sub-manifold $W^s(p)$
there exists a contraction $g: W^s(p)\to W^s(p)$ with fixed point $p_0$
and an injective equivariant immersion $J: W^s(p)\to M$, 
such that $J(p_0)=p$ and $DJ(p_0): T_{p_0}(W^s(p))\to T_p(M)$
is an isomorphism; the image of $J$ defines an immersed submanifold 
$W^s(p)\subset M$ called a {\it stable manifold} of $\varphi$ at $p$. 
Clearly, $dim~(W^s(p))=dim~(V^s)$. 
  We say that $\varphi:M\to M$ is an
{\it Anosov diffeomorphism} ([Anosov 1967]  \cite{Ano1}) if the following condition is
satisfied: there exists a splitting of the tangent bundle
$T(M)$ into a continuous Whitney sum $T(M)=E^s+E^u$, 
invariant under $D\varphi: T(M)\to T(M)$, so that $D\varphi: E^s\to E^s$
is contracting and $D\varphi: E^u\to E^u$ is expanding
%*********************************************************************
\footnote{
It follows from definition, that the Anosov diffeomorphism
imposes a restriction on topology of manifold $M$, in the sense
that not each manifold can support such a diffeomorphism; 
however, if one Anosov diffeomorphism exists on $M$, there 
are infinitely many (conjugacy classes of) such diffeomorphisms
on $M$. It is an open problem of S.~Smale, which $M$ can carry
an Anosov diffeomorphism; so far,  it is proved that    
 the hyperbolic  diffeomorphisms of $m$-dimensional tori  and 
certain automorphisms of the nilmanifolds are Anosov's [Smale 1967]  \cite{Sma1}. 
It is worth mentioning, that on each surface of genus $g\ge 1$ 
there exists a rich family of the so-called pseudo-Anosov diffeomorphisms [Thurston 1988]  \cite{Thu1},
to which our theory  fully applies.  
}
%***********************************************************************
.  Let  $p$ be a fixed point
of the Anosov diffeomorphism $\varphi:M\to M$ and  $W^s(p)$
its stable manifold. Since $W^s(p)$ cannot have self-intersections 
 or limit compacta,   $W^s(p)\to M$ is a dense immersion, i.e. 
 the closure of $W^s(p)$ is the entire $M$.  Moreover, if $q$
is a periodic point of $\varphi$ of period $n$, then  $W^s(q)$ is
a translate of $W^s(p)$, i.e. locally they look like two parallel lines. 
Consider a foliation ${\cal F}$ of $M$, whose leaves are the
translates of $W^s(p)$; the ${\cal F}$ is a continuous foliation
([Smale 1967]  \cite{Sma1}, p.760), which is invariant under the action of diffeomorphism $\varphi$ on its
leaves, i.e. $\varphi$ moves leaves of ${\cal F}$ to the leaves of ${\cal F}$.
The holonomy of ${\cal F}$ preserves the Lebesgue measure and, therefore,
${\cal F}$ is a measured foliation; we shall call it an {\it invariant measured
foliation} and  denote by ${\cal F}_{\varphi}$.
The $AF$-algebra of foliation ${\cal F}$  is called {\it fundamental}, when
${\cal F}={\cal F}_{\varphi}$, where  $\varphi$  is  an Anosov diffeomorphism;
the following is a basic property of such algebras. 
%************************************************************************
\begin{lem}\label{lm2}
Any  fundamental $AF$-algebra is stationary.  
\end{lem}
%***********************************************************************

%**************************************************************************
\section{Proofs}
%***************************************************************************
%**************************************************************************
\subsection{Proof of lemma \ref{lm1}}
%***************************************************************************
Let ${\cal F}'$ be a measured foliation on $M$, given by a closed form
$\omega'\in H^k(M; {\Bbb R})$; let ${\cal F}$ be the measured foliation on $N$,
induced by a submersion $f: N\to M$.  Roughly speaking, we have to prove,
that diagram in Fig.1 is commutative; the proof amounts to the fact, that 
the periods of form $\omega'$ are contained among the periods of form
$\omega\in H^k(N; {\Bbb R})$ corresponding to the foliation ${\cal F}$.  
%*******************************************************************
\begin{figure}[here]
%*******************************************************************
\begin{picture}(300,110)(-100,-5)
\put(20,70){\vector(0,-1){35}}
\put(130,70){\vector(0,-1){35}}
\put(45,23){\vector(1,0){60}}
\put(45,83){\vector(1,0){60}}
\put(15,20){${\Bbb A}_{\cal F}$}
\put(128,20){${\Bbb A}_{\cal F'}$}
\put(17,80){${\cal F}$}
\put(125,80){${\cal F}'$}
\put(60,30){\sf stable}
\put(45,10){\sf homomorphism}
\put(54,90){\sf induction}
\end{picture}
%******************************************************************
\caption{Functor  $F: {\cal M}_0{\cal F}ol\longrightarrow  {\cal R}ng$.}
\end{figure}
%*******************************************************************
The map $f$ defines a homomorphism $f_*: H_k(N)\to H_k(M)$ of the $k$-th
homology groups; let $\{e_i\}$ and $\{e_i'\}$ be a basis in $H_k(N)$ and
$H_k(M)$, respectively. Since $H_k(M)=H_k(N)~/~ker~(f_*)$, we shall denote 
by $[e_i]:= e_i+ker~(f_*)$ a coset representative of $e_i$;  these can be
identified with the elements $e_i\not\in ker~(f_*)$.  The integral
$\int_{e_i}\omega$ defines a scalar product $H_k(N)\times H^k(N; {\Bbb R})\to {\Bbb R}$,
so that $f_*$ is a linear self-adjoint operator; thus, we can write:
%*********************************************************************
\begin{equation}
\lambda_i'=\int_{e_i'}\omega'=\int_{e_i'}f^*(\omega)=\int_{f_*^{-1}(e_i')}\omega=
\int_{[e_i]}\omega\in P({\cal F}),
\end{equation}
%********************************************************************* 
where $P({\cal F})$ is the Plante group (the group of periods) of foliation ${\cal F}$.
Since $\lambda_i'$ are generators of $P({\cal F}')$, we conclude that 
$P({\cal F}')\subseteq P({\cal F})$.  Note,  that $P({\cal F}')=P({\cal F})$
if and only if  $f_*$ is an isomorphism. 

One can apply a criterion of the stable  homomorphism of $AF$-algebras; namely,
${\Bbb A}_{\cal F}$ and ${\Bbb A}_{\cal F'}$ are stably homomorphic, if and only if,
there exists a positive homomorphism $h: G\to H$ between their dimension
groups $G$ and $H$ [Effros 1981]  \cite{E}, p.15. But $G\cong P({\cal F})$ and $H\cong P({\cal F}')$,
while $h=f_*$.  Thus, ${\Bbb A}_{\cal F}$ and ${\Bbb A}_{\cal F'}$ are stably homomorphic.

The functor $F$ is compatible with the composition; indeed,
let $f:N\to M$ and $f':L\to N$ be submersions.  If ${\cal F}$
is a measured foliation of $M$,  one gets the induced 
foliations ${\cal F}'$ and ${\cal F}''$ on $N$ and $L$, respectively;
these foliations fit the diagram
$(L, {\cal F}'')\buildrel f'\over\longrightarrow (N, {\cal F}')\buildrel f\over\longrightarrow (M, {\cal F})$
and the corresponding Plante groups are included: $P({\cal F}'')\supseteq P({\cal F}')\supseteq P({\cal F})$.
Thus, $F(f'\circ f)=F(f')\circ F(f)$, since the inclusion of the Plante groups corresponds to the
composition of homomorphisms;  lemma \ref{lm1} is proved.
$\square$

%**************************************************************************
\subsection{Proof of lemma \ref{lm2}}
%***************************************************************************
Let $\varphi:M\to M$ be an Anosov diffeomorphism; we proceed by showing, that the
invariant foliation ${\cal F}_{\varphi}$ is given by the form $\omega\in H^k(M; {\Bbb R})$,
which is an eigenvector of the linear map $[\varphi]:  H^k(M; {\Bbb R})\to  H^k(M; {\Bbb R})$
induced by $\varphi$.
Indeed, let $0<c<1$ be the contracting constant of the stable sub-bundle $E^s$ of
diffeomorphism $\varphi$ and $\Omega$ the corresponding volume element; by definition,
$\varphi(\Omega)=c\Omega$. Note, that
$\Omega$ is given by restriction of the form $\omega$ to a $k$-dimensional 
manifold, transverse to the leaves of ${\cal F}_{\varphi}$. The leaves of ${\cal F}_{\varphi}$
are fixed by $\varphi$ and, therefore, $\varphi(\Omega)$ is given by a multiple $c\omega$
of form $\omega$. Since $\omega\in H^k(M; {\Bbb R})$ is a vector, whose
coordinates define ${\cal F}_{\varphi}$ up to a scalar, we conclude, that $[\varphi](\omega)=c\omega$,
i.e. $\omega$ is an eigenvector of the linear map $[\varphi]$.  
Let  $(\lambda_1,\dots,\lambda_n)$ be a basis of the Plante group $P({\cal F}_{\varphi})$,
such that $\lambda_i>0$. Notice, that $\varphi$ acts on $\lambda_i$ as multiplication
by constant $c$; indeed, since $\lambda_i=\int_{\gamma_i}\omega$, we have:
%*************************************************************************
\begin{equation}
\lambda_i'=\int_{\gamma_i}[\varphi](\omega)=\int_{\gamma_i}c\omega=c\int_{\gamma_i}\omega=c\lambda_i,
\end{equation}
%**************************************************************
where $\{\gamma_i\}$ is a basis in $H_k(M)$.  Since $\varphi$ preserves the leaves
of ${\cal F}_{\varphi}$, one concludes that $\lambda_i'\in P({\cal F}_{\varphi})$;
therefore, $\lambda_j'=\sum b_{ij}\lambda_i$ for a non-negative 
integer matrix $B=(b_{ij})$. According to [Bauer 1996]  \cite{Bau1}, the matrix $B$ can be 
written as a finite product:
%**************************************************************
\begin{equation}
B=
\left(\matrix{0 & 1\cr I & b_1}\right)\dots
\left(\matrix{0 & 1\cr I & b_p}\right):=B_1\dots B_p,
\end{equation}
%****************************************************************
where $b_i=(b^{(i)}_1,\dots, b^{(i)}_{n-1})^T$ is a vector of non-negative 
integers and   $I$ the unit matrix.   Let $\lambda=(\lambda_1,\dots,\lambda_n)$.  
Consider a purely periodic Jacobi-Perron continued fraction: 
%**************************************************************
\begin{equation}
\lim_{i\to\infty} 
\overline{B_1\dots B_p}
\left(\matrix{0\cr {\Bbb I}}\right),
\end{equation}
%****************************************************************
where  ${\Bbb I}=(0,\dots, 0, 1)^T$;  by a basic property of such fractions, 
it converges to an eigenvector  $\lambda'=(\lambda_1',\dots,\lambda_n')$ of matrix $B_1\dots B_p$
[Bernstein 1971]  \cite{B},  Ch.3.  But $B_1\dots B_p=B$ and $\lambda$ is an eigenvector of matrix $B$;
therefore, vectors $\lambda$ and $\lambda'$ are collinear.  The collinear vectors  
are known to have the same continued fractions;  thus, we have     
%**************************************************************
\begin{equation}
\left(\matrix{1\cr \theta}\right)=
\lim_{i\to\infty} 
\overline{B_1\dots B_p}
\left(\matrix{0\cr {\Bbb I}}\right),
\end{equation}
%****************************************************************
where  $\theta=(\theta_1,\dots,\theta_{n-1})$ and  $\theta_i=\lambda_{i+1}/\lambda_1$. 
Since vector $(1,\theta)$ unfolds into a periodic Jacobi-Perron  continued fraction,
we conclude, that the $AF$-algebra ${\Bbb A}_{\varphi}$ is stationary.
Lemma \ref{lm2} is proved.   
$\square$

%**************************************************************************
\subsection{Proof of theorem \ref{thm1}}
%***************************************************************************
Let $\psi: N\to N$ and $\varphi: M\to M$ be a pair of Anosov diffeomorphisms;
denote by $(N, {\cal F}_{\psi})$ and $(M, {\cal F}_{\varphi})$ the corresponding
invariant  foliations of manifolds $N$ and $M$,  respectively. 
In view of lemma \ref{lm1}, it is sufficient to prove, that the diagram in
Fig.2 is commutative. We shall split  the proof in a series of lemmas.
%*******************************************************************
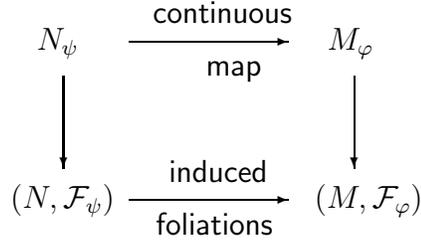
\begin{figure}[here]
%*******************************************************************
\begin{picture}(300,110)(-100,-5)
\put(20,70){\vector(0,-1){35}}
\put(130,70){\vector(0,-1){35}}
\put(45,23){\vector(1,0){60}}
\put(45,83){\vector(1,0){60}}
\put(0,20){$(N, {\cal F}_{\psi})$}
\put(115,20){$(M, {\cal F}_{\varphi})$}
\put(10,80){$N_{\psi}$}
\put(120,80){$M_{\varphi}$}
\put(60,30){\sf induced}
\put(55,10){\sf foliations}
\put(54,90){\sf continuous}
\put(74,70){\sf map}
\end{picture}
%******************************************************************
\caption{Mapping tori and invariant foliations.}
\end{figure}
%*******************************************************************
%************************************************************************
\begin{lem}\label{lm3}
There exists a continuous map $N_{\psi}\to M_{\varphi}$,   whenever
$f\circ\varphi=\psi\circ f$ for a submersion $f: N\to M$.  
\end{lem}
%***********************************************************************
%*******************************************************************
\begin{figure}[here]
%*******************************************************************
\begin{picture}(300,110)(-100,-5)
\put(40,70){\vector(1,-1){35}}
\put(120,70){\vector(-1,-1){35}}
\put(45,83){\vector(1,0){60}}
\put(20,80){$N_{\psi}$}
\put(120,80){$M_{\varphi}$}
\put(75,20){$S^1$}
\put(75,90){$h$}
\put(30,50){$p_{\psi}$}
\put(120,50){$p_{\varphi}$}
\end{picture}
%******************************************************************
\caption{The fiber bundles $N_{\psi}$ and $M_{\varphi}$ over $S^1$.}
\end{figure}
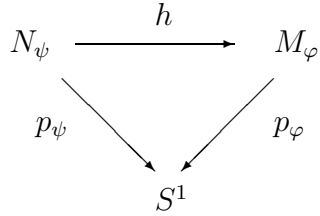
%*******************************************************************
{\it Proof.} (i) Suppose, that $h: N_{\psi}\to M_{\varphi}$ is a continuous
map; let us show, that there exists a submersion $f: N\to M$, such that
$f\circ\varphi=\psi\circ f$.  Both $N_{\psi}$ and $M_{\varphi}$ fiber over
the circle $S^1$ with the projection map  $p_{\psi}$ and $p_{\varphi}$,  respectively; 
therefore, the diagram in Fig.3 is commutative. Let $x\in S^1$; 
since $p_{\psi}^{-1}=N$ and $p_{\varphi}^{-1}=M$, the restriction of $h$
to $x$ defines a submersion $f: N\to M$, i.e. $f=h_x$. Moreover, 
since $\psi$ and $\varphi$ are monodromy maps of  the bundle, it holds:
%*******************************************************************
\begin{equation}
\left\{
\begin{array}{cc}
p_{\psi}^{-1}(x+2\pi)  &= \psi(N),\\
p_{\varphi}^{-1}(x+2\pi)  &= \varphi(M).
\end{array}
\right.
\end{equation}
%*****************************************************************
From the diagram in Fig.3, we get:
$\psi(N)=p_{\psi}^{-1}(x+2\pi)=f^{-1}(p_{\varphi}^{-1}(x+2\pi))=
f^{-1}(\varphi(M))=f^{-1}(\varphi(f(N)))$;  thus, $f\circ\psi=\varphi\circ f$. 
The necessary condition of lemma \ref{lm3} follows.

\bigskip
(ii) Suppose, that $f: N\to M$ is a submersion, such that $f\circ\varphi=\psi\circ f$;
we have to construct a continuous map $h: N_{\psi}\to M_{\varphi}$. Recall,  that
%*******************************************************************
\begin{equation}
\left\{
\begin{array}{cc}
N_{\psi}  &= \{N\times [0,1]~|~(x,0)\sim(\psi(x),1)\}, \\  
M_{\varphi} &= \{M\times [0,1]~|~(y,0)\sim(\varphi(y),1)\}.
\end{array}
\right.
\end{equation}
%***************************************************************** 
We shall identify the points of $N_{\psi}$ and $M_{\varphi}$ using
the substitution $y=f(x)$; it remains to verify, that such an identification
will satisfy the gluing condition $y\sim\varphi(y)$. In view of condition
$f\circ\varphi=\psi\circ f$, we have:
%*****************************************************************
\begin{equation}
y=f(x)\sim f(\psi(x))=\varphi(f(x))=\varphi(y).
\end{equation}
%********************************************************************
Thus, $y\sim\varphi(y)$ and, therefore, the map $h: N_{\psi}\to M_{\varphi}$
is continuous. The sufficient condition of lemma \ref{lm3} is proved.
$\square$

%*******************************************************************
\begin{figure}[here]
%*******************************************************************
\begin{picture}(300,110)(-100,-5)
\put(20,70){\vector(0,-1){35}}
\put(130,70){\vector(0,-1){35}}
\put(45,23){\vector(1,0){60}}
\put(45,83){\vector(1,0){60}}
\put(-10,20){$H^k(M, {\Bbb R})$}
\put(120,20){$H^k(M, {\Bbb R})$}
\put(-10,80){$H^k(N; {\Bbb R})$}
\put(120,80){$H^k(N, {\Bbb R})$}
\put(75,30){$[\varphi]$}
\put(75,90){$[\psi]$}
\put(0,50){$[f]$}
\put(138,50){$[f]$}
\end{picture}
%******************************************************************
\caption{The linear maps $[\psi], [\varphi]$ and $[f]$.}
\end{figure}
%*******************************************************************
%************************************************************************
\begin{lem}\label{lm4}
If a submersion $f: N\to M$ satisfies condition $f\circ\varphi=\psi\circ f$
for the Anosov diffeomorphisms $\psi: N\to N$ and $\varphi: M\to M$,
then the invariant foliations $(N, {\cal F}_{\psi})$ and $(M, {\cal F}_{\varphi})$
are induced by $f$. 
 \end{lem}
%***********************************************************************
{\it Proof.} 
The invariant foliations ${\cal F}_{\psi}$ and ${\cal F}_{\varphi}$
are measured;  we shall denote by $\omega_{\psi}\in H^k(N; {\Bbb R})$
and $\omega_{\varphi}\in H^k(M; {\Bbb R})$ the corresponding cohomology
class,  respectively.  The linear maps on $H^k(N; {\Bbb R})$ and $H^k(M; {\Bbb R})$
induced by $\psi$ and $\varphi$, we shall denote by $[\psi]$ and $[\varphi]$; 
the linear map between $H^k(N; {\Bbb R})$ and $H^k(M; {\Bbb R})$
induced by $f$, we write as $[f]$.  Notice, that $[\psi]$ and $[\varphi]$
are isomorphisms,  while $[f]$ is generally a homomorphism. It was shown
earlier, that $\omega_{\psi}$ and $\omega_{\varphi}$ are eigenvectors of
linear maps $[\psi]$ and $[\varphi]$, respectively; in other words, we have:
%*******************************************************************
\begin{equation}
\left\{
\begin{array}{cc}
\hbox{$[\psi]$}\omega_{\psi}  &=  c_1\omega_{\psi},\\  
\hbox{$[\varphi]$}\omega_{\varphi} &= c_2\omega_{\varphi},
\end{array}
\right.
\end{equation}
%***************************************************************** 
where $0<c_1<1$ and $0<c_2<1$.  Consider a diagram in Fig.4, which involves
the linear maps $[\psi], [\varphi]$ and $[f]$;  the diagram is commutative,
since  condition $f\circ\varphi=\psi\circ f$ implies,  that 
$[\varphi]\circ [f]=[f]\circ [\psi]$. Take the eigenvector $\omega_{\psi}$
and consider its image under the linear map $[\varphi]\circ [f]$:
%*****************************************************************
\begin{equation}
[\varphi]\circ [f](\omega_{\psi})=[f]\circ [\psi](\omega_{\psi})=
[f](c_1\omega_{\psi})=c_1\left([f](\omega_{\psi})\right).
\end{equation}
%*************************************************************
Therefore, vector $[f](\omega_{\psi})$ is an eigenvector of the
linear map $[\varphi]$; let  compare it with the eigenvector $\omega_{\varphi}$:
%*******************************************************************
\begin{equation}
\left\{
\begin{array}{cc}
\hbox{$[\varphi]$}\left(\hbox{$[f]$}(\omega_{\psi})\right)
 &=
 c_1\left(\hbox{$[f]$}(\omega_{\psi})\right),\\ 
\hbox{$[\varphi]$}\omega_{\varphi}  &=  c_2\omega_{\varphi}.
\end{array}
\right.
\end{equation}
%***************************************************************** 
We conclude, therefore, that $\omega_{\varphi}$ and $[f](\omega_{\psi})$
are collinear vectors,  such that $c_1^m=c_2^n$ for some integers $m,n>0$;
a scaling gives us $[f](\omega_{\psi})=\omega_{\varphi}$. 
The latter is an analytic formula, which says that the submersion
$f:N\to M$ induces the foliation $(N, {\cal F}_{\psi})$  from   the
foliation $(M, {\cal F}_{\varphi})$.  Lemma \ref{lm4} is proved.
$\square$

\bigskip
To finish our proof of theorem \ref{thm1}, let $N_{\psi}\to M_{\varphi}$
be a continuous map;  by lemma \ref{lm3}, there exists a submersion
$f: N\to M$, such that $f\circ\varphi=\psi\circ f$. Lemma \ref{lm4}
says  that in this case the invariant measured foliations $(N, {\cal F}_{\psi})$
and $(M, {\cal F}_{\varphi})$ are induced. On the other hand, from lemma \ref{lm2}
we know, that the Jacobi-Perron continued fraction connected to foliations ${\cal F}_{\psi}$
and ${\cal F}_{\varphi}$ is  periodic and, hence,  convergent [Bernstein 1971]  \cite{B};
therefore, one can apply lemma \ref{lm1},   which says  that the $AF$-algebra
${\Bbb A}_{\psi}$ is stably homomorphic to the $AF$-algebra ${\Bbb A}_{\varphi}$.  The latter are,
 by definition, the fundamental $AF$-algebras of the Anosov diffeomorphisms $\psi$ and 
$\varphi$, respectively. Theorem \ref{thm1} is proved.
$\square$

%**************************************************************************
\section{Applications of theorem \ref{thm1}}
%***************************************************************************
%**************************************************************************
\subsection{Galois group of the fundamental $AF$-algebra}
%***************************************************************************
Let ${\Bbb A}_{\psi}$ be  a  fundamental $AF$-algebra and $B$ its primitive incidence matrix,
i.e.  $B$ is not a power of some positive integer matrix.  Suppose that the characteristic polynomial of $B$ is irreducible
and let  $K_{\psi}$ be its  splitting field; then  $K_{\psi}$ is a Galois extension of ${\Bbb Q}$. 
We call $Gal~({\Bbb A}_{\psi}):=Gal~(K_{\psi}|{\Bbb Q})$ the {\it Galois group}  of the fundamental
$AF$-algebra  ${\Bbb A}_{\psi}$; such a group  is determined by the $AF$-algebra  ${\Bbb A}_{\psi}$.   
The second algebraic field is connected to the Perron-Frobenius eigenvalue $\lambda_B$
of the matrix $B$; we shall denote this field ${\Bbb Q}(\lambda_B)$. Note, that
${\Bbb Q}(\lambda_B)\subseteq K_{\psi}$ and ${\Bbb Q}(\lambda_B)$  is not, in general, 
a Galois extension of ${\Bbb Q}$; the obstacle are complex roots of the polynomial 
$char~(B)$ and if there are no such roots  then ${\Bbb Q}(\lambda_B)=K_{\psi}$,
see e.g. [Morandi 1996]  \cite{M}.
There is still a group $Aut~({\Bbb Q}(\lambda_B))$ of automorphisms of ${\Bbb Q}(\lambda_B)$
fixing the field ${\Bbb Q}$ and $Aut~({\Bbb Q}(\lambda_B))\subseteq Gal~(K_{\psi})$
is a subgroup inclusion.  
%*************************************************************************
\begin{lem}\label{lm5}
If $h: {\Bbb A}_{\psi}\to {\Bbb A}_{\varphi}$ is a stable homomorphism, then ${\Bbb Q}(\lambda_{B'})\subseteq K_{\psi}$ is a field inclusion, 
where $B'$ is the matrix associated to $\varphi$.
    \end{lem}
%***********************************************************************
{\it Proof.} Notice, that the non-negative matrix $B$ becomes strictly positive, 
when a proper power of it is taken; we always assume $B$ positive.  Let 
$\lambda=(\lambda_1,\dots,\lambda_n)$ be a basis of the Plante group $P({\cal F}_{\psi})$.
Following the proof of lemma \ref{lm2}, one concludes that $\lambda_i\in K_{\psi}$;
indeed, $\lambda_B\in K_{\psi}$ is the Perron-Frobenius eigenvalue of $B$ , while $\lambda$   
the corresponding eigenvector. The latter can be scaled so, that $\lambda_i\in K_{\psi}$. 
Any stable homomorphism $h: {\Bbb A}_{\psi}\to {\Bbb A}_{\varphi}$ induces a positive
homomorphism of their dimension groups $[h]: G\to H$;  but $G\cong P({\cal F}_{\psi})$
and $H\cong P({\cal F}_{\varphi})$. From inclusion $P({\cal F}_{\varphi})\subseteq P({\cal F}_{\psi})$,
one gets ${\Bbb Q}(\lambda_{B'})\cong P({\cal F}_{\varphi})\otimes {\Bbb Q}\subseteq 
P({\cal F}_{\psi})\otimes {\Bbb Q}\cong {\Bbb Q}(\lambda_B)\subseteq K_{\psi}$ and, therefore, 
${\Bbb Q}(\lambda_{B'})\subseteq K_{\psi}$. Lemma \ref{lm5} follows.
$\square$

%************************************************************************
\begin{cor}\label{cor1}
If $h: {\Bbb A}_{\psi}\to {\Bbb A}_{\varphi}$ is a stable homomorphism, then
$Aut~({\Bbb Q}(\lambda_{B'}))$ (or,  $Gal~({\Bbb A}_{\varphi})$)  is a subgroup 
(or, a  normal subgroup)  of $Gal~({\Bbb A}_{\psi})$.
\end{cor}
%***********************************************************************
{\it Proof.} 
The (Galois) subfields of the Galois field $K_{\psi}$ are bijective with the (normal) subgroups 
of the group $Gal~(K_{\psi})$ [Morandi 1996]  \cite{M}.
$\square$

%**************************************************************************
\subsection{Tight torus bundles}
%***************************************************************************
Let $T^m\cong {\Bbb R}^m/{\Bbb Z}^m$ be an $m$-dimensional torus; let
$\psi_0$ be  a  $m\times m$ integer matrix with $det~(\psi_0)=1$, such that it
is similar to a positive matrix.
The  matrix $\psi_0$ defines  a linear transformation of ${\Bbb R}^m$,
which preserves the lattice $L\cong {\Bbb Z}^m$ of points with integer 
coordinates. There is an induced diffeomorphism $\psi$ of the quotient $T^m\cong {\Bbb R}^m/{\Bbb Z}^m$
onto itself; this diffeomorphism $\psi: T^m\to T^m$ has a fixed point $p$ corresponding 
to the origin of ${\Bbb R}^m$. 
Suppose that $\psi_0$ is hyperbolic, i.e. there are no eigenvalues of $\psi_0$ at the
unit circle; then $p$ is a hyperbolic fixed point of $\psi$ and the stable manifold 
$W^s(p)$ is the image of the corresponding eigenspace of $\psi_0$ under the 
projection ${\Bbb R}^m\to T^m$. If $codim~W^s(p)=1$,  the hyperbolic linear transformation
 $\psi_0$  (and the diffeomorphism $\psi$)  will be  called {\it tight}. 
%*************************************************************************
\begin{lem}\label{lm6}
The tight hyperbolic matrix  $\psi_0$ is similar to the matrix $B$ of the fundamental $AF$-algebra 
${\Bbb A}_{\psi}$.
\end{lem}
%***********************************************************************
{\it Proof.} Since $H_k(T^m; {\Bbb R})\cong {\Bbb R}^{{m!\over k!(m-k)!}}$,
one gets $H_{m-1}(T^m; {\Bbb R})\cong {\Bbb R}^m$; in view of the Poincar\'e 
duality, $H^1(T^m; {\Bbb R})=H_{m-1}(T^m; {\Bbb R})\cong {\Bbb R}^m$. 
Since $codim~W^s(p)=1$, measured foliation ${\cal F}_{\psi}$ is given by a closed
form $\omega_{\psi}\in H^1(T^m; {\Bbb R})$, such that $[\psi]\omega_{\psi}=\lambda_{\psi}\omega_{\psi}$,
where $\lambda_{\psi}$ is the eigenvalue of the linear transformation 
$[\psi]: H^1(T^m; {\Bbb R})\to H^1(T^m; {\Bbb R})$.  It is easy to see that $[\psi]={\psi}_0$, because $H^1(T^m; {\Bbb R})\cong {\Bbb R}^m$
is the universal cover for $T^m$, where  the eigenspace $W^u(p)$ of ${\psi}_0$ is the span of the
eigenform $\omega_{\psi}$. 
On the other hand, from the proof of lemma \ref{lm2} we know that the Plante group $P({\cal F}_{\psi})$ 
is generated by the coordinates of vector $\omega_{\psi}$;  the matrix $B$ is nothing but the matrix
$\psi_0$ written in a new basis of $P({\cal F}_{\psi})$.  
Each change of basis in the ${\Bbb Z}$-module $P({\cal F}_{\psi})$ is given by an integer invertible matrix $S$;
therefore, $B=S^{-1}\psi_0S$. Lemma \ref{lm6} follows.
$\square$

\bigskip\noindent
Let $\psi: T^m\to T^m$ be a hyperbolic diffeomorphism;
the mapping torus $T^m_{\psi}$ will be called a (hyperbolic) {\it torus bundle}
of dimension $m$. Let $k=|Gal~({\Bbb A}_{\psi})|$;  it follows from the Galois theory,
that $1<k\le m!$. Denote by $t_i$ the cardinality of a subgroup $G_i\subseteq Gal~({\Bbb A}_{\psi})$.
%************************************************************************
\begin{cor}\label{cor2}
There are no (non-trivial) continuous map $T^m_{\psi}\to T^{m'}_{\varphi}$,  whenever $t_i'\nmid k$
for all $G_i'\subseteq Gal~({\Bbb A}_{\varphi})$. 
\end{cor}
%***********************************************************************
{\it Proof.} 
If $h: T^m_{\psi}\to T^{m'}_{\varphi}$ was a continuous map to a torus bundle of dimension $m'<m$, then,
by theorem \ref{thm1} and corollary \ref{cor1}, the $Aut~({\Bbb Q}(\lambda_{B'}))$
(or, $Gal~({\Bbb A}_{\varphi})$) were a non-trivial
subgroup (or, normal subgroup)  of the group $Gal~({\Bbb A}_{\psi})$; since $k=|Gal~({\Bbb A}_{\psi})|$, one 
concludes that one of $t_i'$ divides $k$.   This contradicts our assumption.
$\square$
%************************************************************************
\begin{dfn}
The torus bundle $T^m_{\psi}$ is called robust, if there exists $m'<m$, such that
no continuous map $T_{\psi}^m\to T_{\varphi}^{m'}$ exists.  
\end{dfn}
%***********************************************************************
Are there robust bundles? It is shown in this section, that for $m=2,3$ and $4$
there are infinitely many robust torus bundles; a reasonable guess is 
that it is true in any dimension. 
%**************************************************************************
\subsubsection{Case $m=2$}
%***************************************************************************
This case is trivial; $\psi_0$ is a hyperbolic matrix and always tight. 
The $char~(\psi_0)=char~(B)$ is an irreducible quadratic polynomial with two
real roots; $Gal~({\Bbb A}_{\psi})\cong {\Bbb Z}_2$ and, therefore, 
$|Gal~({\Bbb A}_{\psi})|=2$.  Formally, $T^2_{\psi}$ is robust, since no torus 
bundle of a smaller dimension is defined. 
%**************************************************************************
\subsubsection{Case $m=3$}
%***************************************************************************
The $\psi_0$ is hyperbolic; it is always tight, since one root of $char~(\psi_0)$
is real and isolated inside or outside the unit circle.  
%************************************************************************
\begin{cor}\label{cor3}
Let 
%***************************************************************
\begin{equation}
\psi_0(b,c)=\left(\matrix{-b & 1 & 0\cr
                  -c & 0 & 1\cr
                  -1 & 0 & 0}\right)
\end{equation}
%**************************************************************
be such, that $char~(\psi_0(b,c))=x^3+bx^2+cx+1$ is irreducible and $-4b^3+b^2c^2+18bc-4c^3-27$ is 
the square of an integer; then $T^3_{\psi}$ admits no continuous map to any $T_{\varphi}^2$.  
\end{cor}
%***********************************************************************
{\it Proof.} The $char~(\psi_0(b,c))=x^3+bx^2+cx+1$ and the discriminant
 $D=-4b^3+b^2c^2+18bc-4c^3-27$. By Theorem 13.1 [Morandi 1996]  \cite{M}, $Gal~({\Bbb A}_{\psi})\cong {\Bbb Z}_3$
and, therefore, $k=|Gal~({\Bbb A}_{\psi})|=3$. 
For $m'=2$, it was shown that $Gal~({\Bbb A}_{\varphi})\cong {\Bbb Z}_2$ and,
therefore, $t_1'=2$. Since  $2\nmid 3$, corollary \ref{cor2} says that
no continuous map $T_{\psi}^3\to T_{\varphi}^2$ can be constructed.
$\square$

\bigskip\noindent
{\bf Example 1.}
There are infinitely many matrices $\psi_0(b,c)$ satisfying the assumptions of corollary
\ref{cor3}; below are a few numerical examples of robust bundles:
%************************************************************************************
$$   
\left(\matrix{0 & 1 & 0\cr
                  3 & 0 & 1\cr
                  -1 & 0 & 0}\right),
\quad
\left(\matrix{1 & 1 & 0\cr
                  2 & 0 & 1\cr
                  -1 & 0 & 0}\right),
\quad
\left(\matrix{2 & 1 & 0\cr
                  1 & 0 & 1\cr
                  -1 & 0 & 0}\right),
\quad
\left(\matrix{3 & 1 & 0\cr
                  0 & 0 & 1\cr
                  -1 & 0 & 0}\right).
$$
%************************************************************************************
Notice, that the above matrices are not pairwise similar; it can be gleaned from their
traces. Thus, they represent topologically distinct torus bundles. 
%**************************************************************************
\subsubsection{Case $m=4$}
%***************************************************************************
Let $p(x)=x^4+ax^3+bx^2+cx+d$ be a quartic.  Consider the associated cubic polynomial
$r(x)=x^3-bx^2+(ac-4d)x+4bd-a^2d-c^2$;  denote by $L$ the splitting field of $r(x)$. 
%************************************************************************
\begin{cor}\label{cor4}
Let 
%***************************************************************
\begin{equation}
\psi_0(a,b,c)=\left(\matrix{-a & 1 & 0 & 0\cr
                  -b & 0 & 1 & 0\cr
                  -c & 0 & 0 & 1\cr
                  -1 & 0 & 0 & 0}\right)
\end{equation}
%**************************************************************
be tight and such, that $char~(\psi_0(a,b,c))=x^4+ax^3+bx^2+cx+1$ is irreducible and 
one of the following holds: (i) $L={\Bbb Q}$;
(ii) $r(x)$ has a unique root $t\in {\Bbb Q}$ and $h(x)=(x^2-tx+1)[x^2+ax+(b-t)]$
splits over $L$; (iii) $r(x)$ has a unique root $t\in {\Bbb Q}$ and $h(x)$ does not split over $L$.
Then $T_{\psi}^4$ admits no continuous map to any $T_{\varphi}^3$  with $D>0$.  
\end{cor}
%***********************************************************************
{\it Proof.}
According to Theorem 13.4 [Morandi 1996]  \cite{M}, $Gal~({\Bbb A}_{\psi})\cong {\Bbb Z}_2\oplus {\Bbb Z}_2$
in case (i); $Gal~({\Bbb A}_{\psi})\cong {\Bbb Z}_4$ in case (ii); and $Gal~({\Bbb A}_{\psi})\cong D_4$
(the dihedral group) in case (iii). Therefore, $k=|{\Bbb Z}_2\oplus {\Bbb Z}_2|=|{\Bbb Z}_4|=4$
or $k=|D_4|=8$. On the other hand, for $m'=3$ with $D>0$ (all roots are real), we have
$t_1'=|{\Bbb Z}_3|=3$ and $t_2'=|S_3|=6$. Since $3; 6\nmid 4;8$, corollary \ref{cor2} says
that continuous map $T_{\psi}^4\to T_{\varphi}^3$ is impossible.
$\square$

\bigskip\noindent
{\bf Example 2.}
There are infinitely many matrices $\psi_0$,  which satisfy the assumption of 
corollary \ref{cor4};  indeed, consider a family
%***************************************************************
\begin{equation}
\psi_0(a,c)=\left(\matrix{-2a & 1 & 0 & 0\cr
                  -a^2-c^2 & 0 & 1 & 0\cr
                  -2c & 0 & 0 & 1\cr
                  -1 & 0 & 0 & 0}\right),
 \end{equation}
%**************************************************************
where $a, c\in {\Bbb Z}$. The associated cubic  becomes $r(x)=x[x^2-(a^2+c^2)x+4(ac-1)]$,
so that $t=0$ is a rational root; then $h(x)=(x^2+1)[x^2+2ax+a^2+c^2]$. 
The matrix $\psi_0(a,c)$ satisfies one of the conditions (i)-(iii) of 
corollary \ref{cor4} for each $a,c\in {\Bbb Z}$; it remains to 
eliminate the (non-generic) matrices, which are not tight or irreducible.
Thus, $\psi_0(a,c)$ defines  a family of topologically distinct robust bundles.

\bigskip\noindent
%******************************************************************    
{\sf Acknowledgment.} 
I am grateful to the referee for  useful comments.

%**************************************************************************

%**********************************************************

\vskip1cm

\textsc{The Fields Institute for Research in  Mathematical Sciences, Toronto, ON, Canada,  
E-mail:} {\sf igor.v.nikolaev@gmail.com}

%\smallskip
%{\it Current address: Department of Mathematics, 
% University of Sherbrooke, {\sf Igor.Nikolaev@USherbrooke.ca}}  

\end{document}